\documentclass{amsproc}

\usepackage{amsmath, amsthm}%
\usepackage{amsfonts}%
\usepackage{amssymb}%
\usepackage{graphicx}
\usepackage{amscd, amstext, hyperref}
\usepackage{pictexwd}
\usepackage{dcpic}
\usepackage[mathscr]{eucal}
\usepackage{multicol}
\usepackage{multirow}

\usepackage[T1]{fontenc}

\newtheorem{theorem}{Theorem}[section]

\newtheorem{proposition}[theorem]{Proposition}
\newtheorem{corollary}[theorem]{Corollary}

\theoremstyle{definition}
\newtheorem{definition}[theorem]{Definition}
\newtheorem{example}[theorem]{Example}

\theoremstyle{remark}
\newtheorem{remark}[theorem]{Remark}

\numberwithin{equation}{section}

\newcommand{\Z}{\mathbb{Z}}
\newcommand{\Q}{\mathbb{Q}}

\newcommand{\C}{\mathbb{C}}
\newcommand{\F}{\mathbb{F}}

\DeclareMathOperator{\Nm}{Nm}

\newcommand{\Fsf}{\mathsf F}
\newcommand{\Lsf}{\mathsf L}

\newcommand{\dt}{d^T}

\newcommand{\frakp}{\mathfrak{p}}
\newcommand{\legen}[2]{\biggl(\displaystyle{\frac{#1}{#2}}\biggr)}

\newcommand{\psmod}[1]{~(\textup{\text{mod}}~{#1})}

\newcommand{\qq}{q^\times}

\begin{document}

\title{Counting points with Berglund--H\"{u}bsch--Krawitz mirror symmetry}

\author{Ursula Whitcher}
\address{Mathematical Reviews, 416 Fourth St, Ann Arbor, Michigan 48103}
\thanks{The author thanks the Crossing the Walls in Enumerative Geometry conference and the Fields Institute for their hospitality, and the anonymous referee for providing detailed and helpful comments. Edgar Costa and Amanda Francis shared useful and illuminating code. Conversations with Tyler Kelly, Adriana Salerno, and John Voight were, as ever, both delightful and enlightening.}




\begin{abstract}
We give an expository discussion of recent work using Berglund--H\"{u}bsch--Krawitz mirror symmetry to describe the structure of point counts on algebraic varieties over finite fields.
\end{abstract}

\maketitle

\section{Arithmetic mirror symmetry}\label{S:intro}

Can one exploit intuition and constructions from mirror symmetry to prove theorems about the arithmetic of varieties? In this expository paper, we describe work building on observations made by Candelas, de la Ossa, and Rodriguez Villegas around the turn of the millennium. We focus on a particular mirror symmetry construction known as Berglund--H\"{u}bsch--Krawitz mirror symmetry. Our primary aim is to give members of the mirror symmetry research community a taste of the interactions of this subject with number theory. Our approach is generally low-tech and focused on specific examples. We hope that more arithmetically minded readers will appreciate the concrete nature of the constructions, some of which involve varieties of arbitrarily high dimension, and will be inspired to seek unifying principles.

Let us begin by considering a classical example, the Legendre family of elliptic curves, given by the equation:

\[X_{\psi}\colon y^2=x(x-1)(x-\psi).\]

\noindent If we view $X_{\psi}$ as a family of curves defined over the complex numbers, we see that each smooth member of the family admits a holomorphic form which is unique up to scaling.  By choosing a connection on the family, we may fix the choice of scaling and thus the holomorphic form $\omega_\psi$ in a consistent way.  The holomorphic form $\omega_\psi$ satisfies a differential equation depending on the parameter $\psi$, the \emph{Picard--Fuchs differential equation}, which characterizes the changes in the complex structure of the elliptic curves as we move through the family.  

Alternatively, we may view $X_{\psi}$ as a family of elliptic curves defined over a finite field $\mathbb{F}_p$.  In this case, we may consider a different invariant, the \emph{trace of Frobenius} $a_p$, which controls the number of points on $X_\psi$:

\[a_p = 1 + p - \#X_\psi(\mathbb{F}_p).\]

In 1958, Igusa showed that the complex and finite-field aspects of the Legendre family are intimately related.  The Picard--Fuchs equation for the holomorphic form $\omega_\psi$ is \emph{hypergeometric}, with solution 
$$
{}_2F_1\left(\textstyle{\frac{1}{2}, \frac{1}{2}}; 1 \mid \psi \right)=\sum_{n=1}^{\infty} \frac{(\textstyle{\frac{1}{2}})_n^2}{(n!)^2}\psi^n.
$$
Meanwhile, the trace of Frobenius satisfies a truncated hypergeometric formula, with the same coefficients:
$$
a_p \equiv (-1)^{\frac{p-1}{2}} \sum_{n=1}^{(p-1)/2} \frac{(\textstyle{\frac{1}{2}})_n^2}{(n!)^2}\psi^n \pmod{p}.
$$

\noindent Here the generalized hypergeometric equation is given by

\begin{equation}
F(\pmb{\alpha};\pmb{\beta} \mid z) = \sum_{k=0}^{\infty}\frac{(\alpha_1)_k\cdots(\alpha_n)_k}{(\beta_1)_k\cdots(\beta_m)_k}z^k \in \Q[[z]],\label{E:hyperg}
\end{equation}
\noindent where $(x)_k$ is the rising factorial
\[(x)_k = x(x+1)\cdots(x+k-1)=\frac{\Gamma(x+k)}{\Gamma(x)}.\] 

Igusa's result inspired extensive work in number theory highlighting relationships between cohomology, point counting, and hypergeometric structures. We shall focus on the relationship between Igusa's observation and mirror symmetry.

Viewed as a complex manifold, an elliptic curve has a nonvanishing holomorphic form that is unique up to scaling.  Simply connected manifolds with this property are known as \emph{Calabi-Yau manifolds}.  Due to their geometric structure, they play a key role in formulations of string theory.  One may ask whether arithmetic properties of elliptic curves like the ones Igusa observed extend to Calabi-Yau manifolds in higher dimensions.  

Let us consider the diagonal pencil of Calabi-Yau $n-1$-folds in $\mathbb{P}^n$ given by

\[X_{\psi}\colon x_0^{n+1} + \dots + x_n^{n+1} - (n+1) \psi x_0 \cdots x_n =0.\]

\noindent 
Dwork studied this pencil in the case $n=3$ in \cite{dworkpadic}; he described a relationship between the Picard-Fuchs equation and the family's arithmetic properties, which he used to illuminate the structure of his theory of $p$-adic cohomology. Thus, the diagonal pencil in $\mathbb{P}^n$ is usually called the Dwork pencil in the number theory literature, while in the physics literature it is generally referred to as the Fermat pencil. In an ecumenical spirit, and following chronological order, we will refer to it as the Fermat-Dwork pencil.

When $n=4$, the Picard--Fuchs equation for the holomorphic form on $X_{\psi}$ has hypergeometric solution
${}_4F_3\left(\frac{1}{5}, \frac{1}{5}, \frac{1}{5}, \frac{1}{5}; \frac{4}{5}, \frac{3}{5}, \frac{2}{5} \mid \psi^5\right)$.
The hypergeometric structure of this Picard-Fuchs equation is used in the proof of the quintic mirror theorem (see \cite[Chapter 11]{CK} for an expository treatment).  Candelas--de la Ossa--Rodriguez Villegas showed in \cite{CORV} that over a field of prime order, $\#X_\psi(\mathbb{F}_p) \pmod{p}$ is given by a truncation of the hypergeometric series.  

Furthermore, Candelas--de la Ossa--Rodriguez Villegas made a direct link between the arithmetic structure of $X_\psi$ and the structure of the Greene-Plesser mirror $Y_\psi$.  Recall that we can organize information about the number of points on a variety over a finite field in a generating function.  

\begin{definition}
Let $X/\F_q$ be an algebraic variety over the finite field of $q=p^s$ elements.  The \emph{zeta function} of $X$ is 

\[Z(X/\F_q,T):=\exp\left(\sum_{s=1}^{\infty}\#X(\F_{q^s})\frac{T^s}{s}\right)\in \Q[[T]].\]
\end{definition}

For a smooth projective hypersurface $X$ in $\mathbb{P}^n$, we have

\begin{equation}
Z(X, T) = \frac{P_X(T)^{(-1)^n}}{(1-T)(1-qT)\cdots (1-q^{n-1}T)},
\end{equation}

\noindent with $P_X(T) \in \Q[T]$.  The degree of $P_X$ is determined by the Betti numbers of $X$.  If $X$ is a Calabi-Yau hypersurface in $\mathbb{P}^n$, $P_X$ has at most one root that is a $p$-adic unit, termed the \emph{unit root}.  The value of this root determines $\#X(\mathbb{F}_{q}) \pmod{q}$.

If $X$ and $Y$ are mirror Calabi-Yau threefolds, we can expect a relationship between the degrees of certain factors of $Z(X/\F_q,T)$ and $Z(Y/\F_q,T)$ due to the interchange of Hodge numbers.  In some cases one can be more specific.  Candelas, de la Ossa, and Rodriguez Villegas showed in \cite{CORV2} that for the Fermat-Dwork quintic pencil $X_\psi$ and the Greene--Plesser mirror $Y_\psi$, $P_{X_{\psi}}$ and $P_{Y_\psi}$ share a common factor of degree $4$.

We will be particularly interested in detailed analyses of K3 surfaces, so let us move one dimension down and note corresponding results for the K3 surface case.  Let $X_\psi$ be the Fermat--Dwork quartic pencil.  In \cite{kadir}, Shabnam Kadir reports a computation by Xenia de la Ossa, building on Dwork's results, that gives the number of points over a field of prime order $p$:

\begin{equation}
P_{X_\psi}(T)=R_\psi(T)Q_\psi^3(T)S_\psi^{6}(T).
\end{equation} 

\noindent Here (with choices of $\pm$ depending on $p$ and $\psi$), $R_\psi(T)=(1\pm pT)(1-a_\psi T+p^2T)$, $Q_\psi(T)=(1\pm pT)(1\pm pT)$, and $S_\psi(T)=(1-pT)(1+pT)\text{ when $p\equiv3\bmod4$}$ and $(1\pm pT)^2$ otherwise.  We'll discuss point counts for this family over general $\mathbb{F}_q$ in Sections~\ref{S:commonfactor} and \ref{S:pointcount}.

Let $Y_\psi$ be the mirror family to quartics in $\mathbb{P}^3$ (constructed using Greene--Plesser and the Fermat pencil).  Then de la Ossa and Kadir showed:

\begin{equation}
Z(Y_\psi/\F_p,T)=\frac{1}{(1-T)(1-pT)^{19}(1-p^2T)R_\psi(T)}.
\end{equation}

\noindent The factor $R_\psi(T)$ corresponds to periods of the holomorphic form and its derivatives, and is invariant under mirror symmetry.  

One can relate $R_\psi(T)$ to classical arithmetic objects using a natural geometric transformation called a Shioda--Inose structure (see \cite{Dolgachev, ES, Naskrecki}).  Each of the mirror quartics $Y_\psi$ admits an involution $\iota$ which acts symplectically, that is, preserves the holomorphic form.  The quotient $Y_\psi/\iota$ is singular.  Resolving the singularities yields a new K3 surface.  For a general Shioda--Inose structure, the new K3 surface is the Kummer surface associated to an abelian surface.  In the case of mirror quartics, the abelian surface is a product of elliptic curves $E_\psi$ and $E_\psi'$ related by a 2-isogeny:
\[Y_\psi \to Y_\psi/\iota \leftarrow \mathrm{Km}(E_\psi \times E_\psi').\]
\noindent Such pairs of elliptic curves are parametrized by the modular curve $X_0(2)/w_2$, where $w_2$ is an Atkin-Lehner map.  Thus, one can study $R_\psi(T)$ using properties of the modular curve.

One immediately asks how far the point-counting results for the Fermat-Dwork quartic and quintic pencils and their mirrors can be generalized.  Daqing Wan considered the Fermat-Dwork pencils $X_\psi$ and their mirrors $Y_\psi$ in any dimension.  In \cite{wan}, he showed that for any $\psi$ where $X_\psi$ is smooth, the unit roots of $X_\psi$ and $Y_\psi$ coincide.  Thus, we obtain a congruence of point counts:
\begin{equation}\label{E:strongmirror}
\#X_\psi(\mathbb{F}_{q}) \equiv \#Y_\psi(\mathbb{F}_{q}) \pmod{q}.
\end{equation}

Kadir studied a two-parameter family of octic Calabi-Yau threefolds and its generalized Greene--Plesser mirror in \cite{kadir} and \cite{kadir2}.  She showed that over fields of prime order, the zeta functions of the octic threefolds and their mirrors have a common factor, using an explicit point count and techniques from toric varieties.  In \cite{kadir2}, Kadir argues that a similar computation will show in general that if $X_{\psi_1,\dots,\psi_r}$ is a generalized Fermat pencil of Calabi-Yau varieties in a Gorenstein Fano weighted projective space and $Y_{\psi_1,\dots,\psi_r}$ is the Greene-Plesser mirror, $Z(X_{\psi_1,\dots,\psi_r}/\F_p,T)$ and $Z(Y_{\psi_1,\dots,\psi_r}/\F_p,T)$ share a common factor. 

Building on such analysis, Kloosterman studied deformations of generalized Fermat hypersurfaces in weighted projective spaces by monomials other than $x_0 \cdots x_n$ in \cite{kloostermanFermat}. He showed that the corresponding Picard--Fuchs equations are of generalized hypergeometric form, and used this result to analyze the factorization structure of the zeta function. More general deformations of diagonal hypersurfaces are also important for purely arithmetic reasons; for example, in \cite{PT} Pancratz and Tuitman describe an algorithm to compute the zeta function of a projective hypersurface using a $p$-adic version of the Picard--Fuchs differential equation for a one-parameter deformation of a diagonal hypersurface.

The evidence for arithmetic mirror symmetry phenomena for Greene--Plesser mirror symmetry raises the question of whether other mirror constructions also have arithmetic implications.  In \cite{MW}, Magyar and the present author described experimental point-counting results for Calabi-Yau hypersurfaces in Fano toric varieties and gave a conjectural characterization of pencils satisfying a congruence like Equation~\ref{E:strongmirror}.  The authors of \cite{HLYY} use Hasse--Witt matrices to extract information about the unit roots of toric Calabi-Yau hypersurfaces.

In the following sections, we discuss arithmetic implications of Berglund--H\"{u}bsch--Krawitz (BHK) mirror symmetry.  We begin by reviewing the BHK mirror construction in Section~\ref{S:BHK}.  In Section~\ref{S:affine}, we describe a combinatorial trick from \cite{perunicic} for counting points on the affine variety determined by an invertible polynomial. Though elementary, the method draws on aspects of the BHK construction in a way that hints at a deeper cohomological interpretation. In Section~\ref{S:commonfactor}, we use intuition from mirror symmetry to identify common factors in zeta functions of different Calabi-Yau pencils.  Our exposition follows \cite{zeta}, though we also discuss work of Kloosterman in \cite{kloosterman}.  We focus on a specific set of K3 surface examples in Section~\ref{S:pointcount}; we exploit the hypergeometric structure of Picard--Fuchs equations and point counts to describe an explicit motivic deconstruction for these families, following \cite{hypergeometric}.

\section{Berglund-H\"{u}bsch-Krawitz mirrors}\label{S:BHK}

Let us recall the procedure for the Berglund-H\"{u}bsch-Krawitz (BHK) mirror symmetry construction.  Consider a polynomial $F_A$ that is the sum of $n+1$ monomials in $n+1$ variables

\begin{equation}
F_A = \sum_{i=0}^n \prod_{j=0}^n x_j^{a_{ij}}.
\end{equation}

\noindent We view $F_A$ as determined by an integer matrix $A = (a_{ij})$, so each row of the matrix corresponds to a monomial.

\begin{definition}
We say $F_A$ is \emph{invertible} if the matrix $A$ is invertible, there exist positive integers called \emph{weights} $r_j$ so that $d : = \sum_{j=0}^n r_j a_{ij}$ is the same constant for all $i$, and the polynomial $F_A$ has exactly one critical point, namely at the origin. 
\end{definition}

\begin{definition}
We say an invertible polynomial $F_A$ satisfies the \emph{Calabi-Yau condition} if $d= \sum_{j=0}^n r_j$.
\end{definition}

If a polynomial is invertible and the Calabi-Yau condition is satisfied, the weights determine a weighted projective space $\mathbb{WP}^n(r_0,\ldots, r_n)$ and $F_A$ determines a Calabi-Yau hypersurface $X_A$ in this weighted projective space. Alternatively, one may view $F_A$ as a function from $\mathbb{C}^{n+1}$ to $\mathbb{C}$; this point of view yields a \emph{Landau--Ginzburg model}.

Kreuzer and Skarke proved in \cite{KS} that any invertible polynomial $F_A$ can be written as a sum of invertible potentials, each of which must be of one of the three \emph{atomic types}:

\begin{equation}\begin{aligned}
W_{\text{Fermat}} &:= x^a, \\
W_{\text{loop}} &:= x_1^{a_1}x_2 + x_2^{a_2}x_3 + \ldots +x_{m-1}^{a_{m-1}}x_m + x_m^{a_m}x_1, \text{ and } \\
W_{\text{chain}} &:=x_1^{a_1}x_2 + x_2^{a_2}x_3 + \ldots x_{m-1}^{a_{m-1}}x_m + x_m^{a_m}.
\end{aligned}\end{equation}

We will be interested in three subgroups of the torus $(\C^*)^{n+1}$ related to $F_A$. Let $\mathrm{Aut}(A)$ be the diagonal symmetries of $F_A$, that is, those $(\lambda_0, \dots, \lambda_n) \in (\C^*)^{n+1}$ such that $F_A(\lambda_0 x_0, \dots, \lambda_n x_n) = F_A (x_0, \dots, x_n)$ for all $(x_0, \dots, x_n) \in \C^{n+1}$. Note that $\mathrm{Aut}(F_A)$ is a finite abelian group, and the coordinates of each element of $\mathrm{Aut}(F_A)$ are roots of unity. One may compute $\mathrm{Aut}(F_A)$ by using the fact that it is generated by the columns of the matrix $A^{-1}$: if $r_{ij}$ are the coordinates of column $\rho_j$, then $(e^{2\pi i r_{0j}}, \dots, e^{2\pi i r_{nj}})$ is the corresponding element of $\mathrm{Aut}(F_A)$. Let $SL(F_A) \subset (\C^*)^{n+1}$ be the diagonal symmetries of $F_A$ with the property that the product of coordinates $\lambda_0\cdots\lambda_n$ is $1$. Let $J(F_A)$ be the trivial diagonal symmetries, that is, the subgroup of $SL(F_A)$ that acts trivially on $\mathbb{WP}^n(r_0,\ldots, r_n)$.  Then $SL(F_A)/J(F_A)$ acts nontrivially and \emph{symplectically} on $X_A$; that is, it fixes the holomorphic $n-1$-form.

Now, suppose we start with a Calabi-Yau manifold $X_A$ corresponding to an invertible matrix $A$. To construct the Berglund--H\"{u}bsch--Krawitz (BHK) mirror of $X_A$, we take the transpose matrix $A^T$.  This matrix is also invertible; we refer to its weights as \emph{dual weights}.  Consider the polynomial $F_{A^T}$. Let $\widetilde{G^T} = SL(F_{A^T})/J(F_{A^T})$. We obtain a dual orbifold $X_{A^T} / \widetilde{G^T}$ as the mirror of $X_A$. In general, for any group $H$ such that $J(F_A) \subset H \subset SL(F_A)$, one may define the Berglund-H\"{u}bsch-Krawitz mirror of the orbifold $X_A/\widetilde{H}$.  Under this more general construction, BHK duality is a true duality: the mirror of the mirror yields the original orbifold.  Furthermore, BHK duality for the polynomial $F_{A}$ extends naturally to the pencil of hypersurfaces described by
\[F_{A}-(d^T)\psi x_0 \cdots x_n,\]
where $d^T=\sum q_i$ is the sum of the dual weights.

Berglund--H\"{u}bsch--Krawitz mirror symmetry involves correspondences of algebras, as well as spaces. Fan, Jarvis, and Ruan constructed a Gromov--Witten style theory, FJRW theory, for Landau--Ginzburg models; see \cite{FJR} for a detailed discussion. This theory involves Frobenius algebras constructed from the elements of $G$. Chiodo and Ruan proved in \cite{CR} that there is an isomorphism between the FJRW theory of $X/G$ and the Chen-Ruan orbifold cohomology of $X_A/\widetilde{G}$. When $X_A/\widetilde{G}$ and $X_{A^T}/\widetilde{G^T}$ admit crepant resolutions, this isomorphism can be used to recover the classical mirror relationship between Hodge diamonds.

\section{An affine point count}\label{S:affine}

In his dissertation, \cite{perunicic}, Peruni\v{c}i\'{c} showed that if $\det A \mid (p-1)$, then there is a formula for the number of points on the affine variety cut out by $F_A$ in $\mathbb{F}_p^{n+1}$ $\pmod{p}$ that depends only on the matrix $A^T$. Though the argument reduces to elementary combinatorics, it highlights objects of interest in BHK mirror symmetry and hints at deeper connections between the arithmetic of $X_A$ and the structure of its mirror. Let us discuss the details, and correct a small error in \cite{perunicic}, as a warmup.

\begin{proposition}[Theorem 4.3.2, \cite{perunicic}]\label{P:affine}
Let $F_A$ be an invertible polynomial, and suppose $\det A \mid (p-1)$. For any rational vector $\vec{\xi} = (\xi_0, \dots, \xi_n)$, set $\mathrm{age}(\vec{\xi}) = \sum \xi_i$ and
\[\nu(\xi) = \frac{(p-1)!}{\prod_{i=0}^n ((p-1)\xi_i)!}.\] 
Let $\Xi$ be the set of positive integer linear combinations of the columns of $(A^T)^{-1}$:

\[\Xi = \{(A^T)^{-1} \vec{v} \mid v_0, \dots, v_n \in \mathbb{Z}^+ \}.\]
Then the number of points $\nu(A)$ on the affine hypersurface $\{ \vec{x} \in \mathbb{F}_p^n \mid F_a(\vec{x}) = 0\}$ satisfies the equivalence
\[\nu(A) \equiv (-1)^n \sum_{\substack{\xi \in \Xi \\ \mathrm{age}(\xi) = 1}} \nu(\xi) \pmod{p}.\]
\end{proposition}

The proof of Proposition~\ref{P:affine} uses the fact that $\nu \equiv \sum_{\vec{x} \in \mathbb{F}_{p}^n} (1-F_{A}(\vec{x})^{p-1}) \pmod{p}$, together with the multinomial theorem. Note that the Calabi-Yau condition is not required: though we will focus on that case here, Proposition~\ref{P:affine} applies to more general invertible polynomials.

\begin{example}\label{E:ec}
Let $F_A = x_0^2 x_1 + x_1^2 x_2 + x_2^3$. The weights of $A$ are $(1,1,1)$, and the weights of $A^T$ are $(2,1,1)$. We have $\det A = 12$; let us take $p=73$.
The set $\Xi$ contains one element, the vector $(\frac{1}{2}, \frac{1}{4}, \frac{1}{4})$. We compute:
\begin{align*}
\nu(A) & \equiv (-1)^2 \frac{72!}{(72/2)!(72/4)!(72/4)!} \pmod{73}\\
& \equiv 67 \pmod{73}.
\end{align*}
One can check that $\nu(A) = 5761$, which is indeed equivalent to $67 \pmod{73}$.
\end{example}

\begin{example}\label{E:affineK3}
Let $F_A = x_0^2 x_1 + x_1^5 + x_2^5 + x_3^5$. The weights of $A$ are $(2,1,1,1)$, and the weights of $A^T$ are $(5,1,2,2)$. We have $\det A = 250$; let us take $p=251$. The set $\Xi$ contains one element, the vector $(\frac{1}{2}, \frac{1}{10}, \frac{1}{5}, \frac{1}{5})$. We compute:
\begin{align*}
\nu(A) & \equiv (-1)^3 \frac{250!}{(250/2)!(250/10)!(250/5)!(250/5)!} \pmod{251}\\
& \equiv 6 \pmod{251}.
\end{align*}
\end{example}

Since the columns of $(A^T)^{-1}$ generate the group $\mathrm{Aut}(A^T)$, it is natural to try to rephrase Proposition~\ref{P:affine} in terms of elements of $\mathrm{Aut}(A^T)$.

\begin{definition}
Let $g= (e^{2\pi i \xi_0}, \dots, e^{2\pi i \xi_n})$ be a diagonal symmetry written in canonical form, that is, assume $0 \leq \xi_i < 1$ for each $\xi_i$. Then the \emph{age} of $g$ is given by
\[\mathrm{age}(g) = \sum_{i=0}^n \xi_i.\]
\end{definition}

\noindent We shall identify $g$ with the corresponding $\xi$ vector $(\xi_0, \dots, \xi_n)$, when confusion will not arise. 

\begin{remark}
The version of Proposition~\ref{P:affine} given in Peruni\v{c}i\'{c}'s dissertation (\cite[Theorem 4.3.2]{perunicic}) takes a sum over all group elements $g \in \mathrm{Aut}(A^T)$ such that $\mathrm{age}(g) = 1$. However, this statement is in error: such a sum can lead to extraneous terms. For example, if we take $F_A = x_0^2 x_1 + x_1^2 x_2 + x_2^3$ as in Example~\ref{E:ec}, then $\mathrm{Aut}(A^T)$ contains the group element given by $(0, \frac{1}{2}, \frac{1}{2})$, which contributes a gratuitous $\frac{(p-1)!}{(p-1)!^2/4}$ to the sum.
\end{remark}

\begin{definition}
If $g \in \mathrm{Aut}(A)$ acts nontrivially on each of the variables $x_0, \dots, x_n$, that is, the corresponding $\xi$ vector contains no zeros, then we say $g$ is \emph{narrow}. Otherwise, we say $g$ is \emph{broad}.
\end{definition}

If $\xi \in \Xi$, the element $g$ of $\mathrm{Aut}(A)$ determined by $\xi$ must be narrow, because the columns of $(A^T)^{-1}$ are linearly independent. When $A$ (and thus $A^T$) satisfies the Calabi-Yau condition, there is precisely one element of $\Xi$, the vector given by $(\frac{q_0}{d^T}, \dots, \frac{q_n}{d^T})$, where $q_0, \dots, q_n$ are the dual weights. The corresponding group element generates $J(F_{A^T})$, the trivial symmetries of the associated weighted projective space.

Chiodo and Ruan observe in \cite{CR} that, under their Landau--Ginzburg/Calabi--Yau correspondence, each narrow element of $G^T$ yields a cohomology class generated by hyperplanes in $H^{p,q}_{CR}(X_{A^T}/\widetilde{G^T})$. In particular, for any $G$ satisfying $J(F_A) \subset G \subset SL(F_A)$, we have  $J(F_{A^T}) \subset G^T$. Because it has age 1, the generator of $J(F_A)$ given by $(\frac{q_0}{d^T}, \dots, \frac{q_n}{d^T})$ yields an element of $H^{n-1,n-1}_{CR}(X_{A^T}/\widetilde{G^T})$. Under mirror symmetry, we obtain a corresponding element of $H^{n-1,1}_{CR}(X_{A}/\widetilde{G})$. In particular, when $G = J(F_A)$, we may work in ordinary rather than Chen--Ruan cohomology. In this case, Proposition~\ref{P:affine} tells us that for certain primes $p$, we may compute the number of points on the affine cone over $X_A$ $\pmod{p}$, and thus the number of points on $X_A$ $\pmod{p}$, in terms of information associated to a particular class in $H^{n-1,1}(X_{A})$. In the following sections, we shall observe such an association in a different way, by noting the relationship between a particular factor of the zeta function whose root is the unit root and an element of $H^{n-1,1}(X_{A})$ corresponding to a deformation of the holomorphic form. We will also give more general, hypergeometric formulas for many point counts.

In \cite{AP}, Aldi and Peruni\v{c}i\'{c} pursue another strategy for studying the arithmetic and cohomological structure of the mirror correspondence. Borisov gave a vertex algebra formulation of BHK mirror symmetry in \cite{borisovVA} that unifies rings associated to both the A- and B-models in a single algebraic structure. Aldi and Peruni\v{c}i\'{c} realize this vertex algebra construction in the setting of $p$-adic D-modules, and show that it can be made compatible with the Frobenius action. The challenge is then to link this structure to specific arithmetic or geometric predictions.

\section{Common factors}\label{S:commonfactor}

In the early 1990s, Greene, Plesser, and Roan showed in \cite{GPR} that one may construct a mirror family to smooth quintics in $\mathbb{P}^4$ using discrete group quotients of pencils other than the Fermat pencil.  Similarly, the Fermat pencil is not the only highly symmetric pencil one can use to construct the mirror to smooth quartics in $\mathbb{P}^3$.  The alternatives are certain invertible pencils; we list them in Table~\ref{Ta:quartic}.

\begin{table}[h!]
\begin{tabular}{c|c|c}
 Family & Equation	&  $SL(F_A)/J(F_A)$ \\
\hline	\hline
\rule{0pt}{2.5ex}   $\Fsf_4$ & $x_0^4+x_1^4 + x_2^4 + x_3^4 - 4\psi x_0x_1x_2x_3$ 	& $(\Z/4\Z)^2$	 \\
 $\Fsf_2\Lsf_2$ & $x_0^4 + x_1^4 + x_2^3x_3 + x_3^3x_2 - 4\psi x_0x_1x_2x_3$	& $\Z/8\Z$ \\
 $\Fsf_1\Lsf_3$ & $x_0^4 + x_1^3x_2 + x_2^3x_3 + x_3^3x_1 - 4\psi x_0x_1x_2x_3$ & $\Z / 7\Z$  \\
 $\Lsf_2\Lsf_2$ & $x_0^3x_1 + x_1^3x_0 + x_2^3x_3 + x_3^3x_2 - 4\psi x_0x_1x_2x_3$ & $\Z/4\Z \times \Z/2\Z$ \\
 $\Lsf_4$ & $x_0^3x_1 + x_1^3x_2 + x_2^3x_3 + x_3^3 x_0 - 4\psi x_0x_1x_2x_3$ & $\Z / 5\Z$  \\
\end{tabular}
\caption{Symmetric quartic pencils}\label{Ta:quartic}
\end{table}

The zeta functions of these pencils were studied in \cite{zeta}.  For each $\psi$ such that the corresponding K3 surfaces $X_{\diamond, \psi}$ are smooth and nondegenerate, the zeta functions $Z(X_{\diamond,\psi}/\F_q,T)$ share a common factor  $R_{\psi}(T)$ of degree 3 (see \cite[Theorem 5.1.3]{zeta}, and note that we have suppressed the dependence on $q$).  

For $\F_q$ containing sufficiently many roots of unity, it holds that 

\begin{equation}
Z(X_{\diamond,\psi}/\F_q,T)=\frac{1}{(1-T)(1-qT)^{19}(1-q^2T)R_{\psi}(T)}. 
\end{equation}

\noindent We say the zeta functions $Z(X_{\diamond,\psi})$ are \emph{potentially equal}.  Let $Y_\psi$ be the family of mirror quartics constructed by the Greene--Plesser quotient of the Fermat pencil.  Then $Z(X_{\diamond,\psi})$ and $Z(Y_\psi)$ are potentially equal for any $\diamond$.

In this case, we have used the physical intuition that the pencils $X_{\diamond,\psi}$ have the same mirror, in the sense of Greene--Plesser--Roan, to extract arithmetic consequences.  Since all of the pencils listed in Table~\ref{Ta:quartic} are invertible, one may ask whether this phenomenon holds for other cases of the BHK mirror construction: if $X_{A^T}$ and $X_{B^T}$ have common properties, do $X_A$ and $X_B$ share arithmetic properties? The results of \cite{zeta} show that for projective invertible pencils, one can detect common properties using only the matrices $A^T$ and $B^T$:

\begin{theorem}\cite{zeta}\label{T:commonfactor}
Let $X_{A,\psi}$ and $X_{B,\psi}$ be invertible pencils of Calabi-Yau $(n-1)$-folds in $\mathbb{P}^n$.  Suppose $A$ and $B$ have the same dual weights $(q_0, \dots, q_n)$, and let $d^T=q_0+\dots+q_n$. Then for each $\psi \in \mathbb{F}_q$ such that $\gcd(q,(n+1)d^T)=1$ and the fibers $X_{A,\psi}$ and $X_{B,\psi}$ are nondegenerate and smooth, 
the polynomials $P_{X_{A, \psi}}(T)$ and $P_{X_{B,\psi}}(T)$ have a common factor $R_\psi(T) \in \mathbb{Q}[T]$ with 
\[ \deg R_{\psi}(T) \geq D(q_0,\dots,q_n),\]
where $D(q_0,\dots,q_n)$ is the degree of the Picard--Fuchs equation for the holomorphic form.
\end{theorem} 

By results of G\"{a}hrs in \cite{gahrs}, the Picard--Fuchs equation for the holomorphic form depends only on the dual weights, and is hypergeometric. The proof of Theorem~\ref{T:commonfactor} uses this fact together with Dwork's $p$-adic cohomology theory. Though the details of the proof depend on properties of exponential sums, the intuition is that the zeta function can be calculated as the characteristic polynomial of the Frobenius action on $p$-adic cohomology, and that one can identify a subspace of this cohomology space corresponding to the holomorphic form and its derivatives.

One may also use $A^T$ and $B^T$ to detect coincidences of unit roots. This gives less information about the structure of the zeta function, but requires fewer conditions on the field $\mathbb{F}_q$.

\begin{proposition}\cite{zeta}
Let $F_A(x)$ and $F_B(x)$ be invertible polynomials in $n+1$ variables satisfying the Calabi--Yau condition.  Suppose $A^T$ and $B^T$ have the same weights.  Then for all $\psi \in \mathbb{F}_q$ and in all characteristics (including when $p \mid d^T$), either the unit root of $X_{A,\psi}$ is the same as the unit root of $X_{B,\psi}$, or neither variety has a nontrivial unit root. 
\end{proposition}

\begin{corollary}\cite{zeta}
Let $F_A(x)$ and $F_B(x)$ be invertible polynomials in $n+1$ variables satisfying the Calabi--Yau condition.  Suppose $A^T$ and $B^T$ have the same weights.  Then for any fixed $\psi \in \mathbb{F}_q$ and in all characteristics (including $p \mid d^T$) the $\mathbb{F}_q$-rational point counts for fibers $X_{A,\psi}$ and $X_{B,\psi}$ are congruent as follows: 
\[\#X_{A,\psi} \equiv \#X_{B,\psi} \pmod{q}.\]
\end{corollary}

In these cases, the unit root is determined by a formal power series depending on
\[
{}_{D}F_{D-1}\left(\alpha_i; \beta_j \bigm| ({\prod}_i q_i^{-q_i}) \psi^{-\dt}\right)\] 
where the parameters $\alpha_i$ and $\beta_j$ depend only on the dual weights $q_i$. This follows from results of Miyatani in \cite{Miyatani} when $X_{A,\psi}$ is smooth and $\psi \neq 0$, or Adolphson and Sperber in \cite{AS}, in general. Thus, there is a truncated hypergeometric formula for $\#X_\psi(\mathbb{F}_q) \pmod{q}$, generalizing the results of Igusa for the Legendre family.

Kloosterman showed in \cite{kloosterman} that one can extend Theorem~\ref{T:commonfactor} to a broader class of hypersurfaces: he allows for invertible polynomials $A$ and $B$ that do not necessarily satisfy the Calabi-Yau condition, and permits more general one-parameter monomial deformations.  The resulting common factor of $P_{X_{A, \psi}}(T)$ and $P_{X_{B,\psi}}(T)$ may be of larger degree than the common factor identified in Theorem~\ref{T:commonfactor}. 

The idea of Kloosterman's proof is to use a Fermat space $Y_\psi$ covering both $X_{A,\psi}$ and $X_{B,\psi}$, rather than relying on the hypergeometric structure. This construction, which originated in \cite{shioda}, is called a Shioda map. In mirror symmetry, it has been used to study ``multiple mirror'' phenomena, as in \cite{kelly}. Kloosterman studies Shioda maps in more generality, for deformations of invertible polynomials corresponding to hypersurfaces in weighted projective spaces, and shows that when $Y_\psi$ is smooth, the characteristic polynomial of Frobenius acting on a certain subspace of $H^{n-1}(Y_\psi)$ divides the characteristic polynomials of Frobenius on both $H^{n-1}(X_{A,\psi})$ and $H^{n-1}(X_{B,\psi})$. In the projective case, no cancellation occurs, so the result for characteristic polynomials yields a factor of the zeta function.

\section{Point counting and hypergeometric formulas}\label{S:pointcount}

Theorem~\ref{T:commonfactor} shows that, for invertible pencils describing hypersurfaces $X$ in $\mathbb{P}^n$, the piece of middle cohomology corresponding to the holomorphic form and its derivatives corresponds to a factor of $P_X$. This raises the question of whether one can describe other factors of $P_X$ in a similar fashion. For the families of K3 surfaces $X_{\diamond,\psi}$ with $\diamond \in \{\Fsf_4,\Fsf_2\Lsf_2,\Fsf_1\Lsf_3,\Lsf_2\Lsf_2,\Lsf_4\}$, this correspondence is worked out completely in \cite{hypergeometric}.

Let us take the Dwork-Fermat quartic pencil $\Fsf_4$ and the $\Lsf_2\Lsf_2$ family as examples.  The polynomials $P_X$ for several different values of $\psi$ when $p=q=281$ are given in Table~\ref{Ta:F4zeta} and Table~\ref{Ta:L2L2zeta}.  We computed this data using code written by Edgar Costa and described in \cite{CT}.

\begin{table}[ht]
\begin{tabular}{|c|c|}
\hline
$\psi$ & $\Fsf_4$\\
\hline
0& $ ( 1 - 281 T ) ^{ 19 } ( 1 + 462 T + 281 ^ { 2 } T ^ { 2 } ) $ \\ 
1, 53, 228, 280&not smooth\\
2, 106, 175, 279& $ ( 1 - 281 T ) ^{ 3 } ( 1 + 281 T ) ^{ 16 } ( 1 + 238 T + 281 ^ { 2 } T ^ { 2 } ) $ \\ 
3, 122, 159, 278& $ ( 1 - 281 T ) ^{ 19 } ( 1 + 78 T + 281 ^ { 2 } T ^ { 2 } ) $ \\ 
4, 69, 212, 277& $ ( 1 - 281 T ) ^{ 3 } ( 1 + 281 T ) ^{ 16 } ( 1 - 434 T + 281 ^ { 2 } T ^ { 2 } ) $ \\ 
5, 16, 265, 276& $ ( 1 - 281 T ) ^{ 13 } ( 1 + 281 T ) ^{ 6 } ( 1 + 418 T + 281 ^ { 2 } T ^ { 2 } ) $ \\ 
6, 37, 244, 275& $ ( 1 - 281 T ) ^{ 3 } ( 1 + 281 T ) ^{ 16 } ( 1 - 50 T + 281 ^ { 2 } T ^ { 2 } ) $ \\ 
7, 90, 191, 274& $ ( 1 - 281 T ) ^{ 3 } ( 1 + 281 T ) ^{ 16 } ( 1 + 238 T + 281 ^ { 2 } T ^ { 2 } ) $ \\ 
8, 138, 143, 273& $ ( 1 - 281 T ) ^{ 3 } ( 1 + 281 T ) ^{ 16 } ( 1 - 50 T + 281 ^ { 2 } T ^ { 2 } ) $ \\ 
9, 85, 196, 272& $ ( 1 - 281 T ) ^{ 3 } ( 1 + 281 T ) ^{ 16 } ( 1 - 50 T + 281 ^ { 2 } T ^ { 2 } ) $ \\ 
10, 32, 249, 271& $ ( 1 - 281 T ) ^{ 5 } ( 1 + 281 T ) ^{ 16 } $ \\ 
\hline
\end{tabular}
\caption{Examples of $P_X$ for $\Fsf_4$ when $q=281$}\label{Ta:F4zeta}
\end{table}

\begin{table}[ht]
\begin{tabular}{|c|c|}
\hline
$\psi$ & $P_X$\\
\hline
0& $ ( 1 - 281 T ) ^{ 19 } ( 1 + 462 T + 281 ^ { 2 } T ^ { 2 } ) $ \\ 
1, 53, 228, 280&not smooth\\
2, 106, 175, 279& $ ( 1 - 281 T ) ^{ 11 } ( 1 + 281 T ) ^{ 8 } ( 1 + 238 T + 281 ^ { 2 } T ^ { 2 } ) $ \\ 
3, 122, 159, 278& $ ( 1 - 281 T ) ^{ 19 } ( 1 + 78 T + 281 ^ { 2 } T ^ { 2 } ) $ \\ 
4, 69, 212, 277& $ ( 1 - 281 T ) ^{ 15 } ( 1 + 281 T ) ^{ 4 } ( 1 - 434 T + 281 ^ { 2 } T ^ { 2 } ) $ \\ 
5, 16, 265, 276& $ ( 1 - 281 T ) ^{ 13 } ( 1 + 281 T ) ^{ 6 } ( 1 + 418 T + 281 ^ { 2 } T ^ { 2 } ) $ \\ 
6, 37, 244, 275& $ ( 1 - 281 T ) ^{ 15 } ( 1 + 281 T ) ^{ 4 } ( 1 - 50 T + 281 ^ { 2 } T ^ { 2 } ) $ \\ 
7, 90, 191, 274& $ ( 1 - 281 T ) ^{ 11 } ( 1 + 281 T ) ^{ 8 } ( 1 + 238 T + 281 ^ { 2 } T ^ { 2 } ) $ \\ 
8, 138, 143, 273& $ ( 1 - 281 T ) ^{ 15 } ( 1 + 281 T ) ^{ 4 } ( 1 - 50 T + 281 ^ { 2 } T ^ { 2 } ) $ \\ 
9, 85, 196, 272& $ ( 1 - 281 T ) ^{ 15 } ( 1 + 281 T ) ^{ 4 } ( 1 - 50 T + 281 ^ { 2 } T ^ { 2 } ) $ \\ 
10, 32, 249, 271& $ ( 1 - 281 T ) ^{ 17 } ( 1 + 281 T ) ^{ 4 } $ \\ 
\hline
\end{tabular}
\caption{Examples of $P_X$ for $\Lsf_2\Lsf_2$ when $q=281$}\label{Ta:L2L2zeta}
\end{table}

One can compute a Picard--Fuchs equation for any primitive differential form, not just the holomorphic form. We obtain hypergeometric (or trivial) differential equations organized by the group of symmetries. 

The Picard--Fuchs equations for the Dwork-Fermat quartic are given by the following hypergeometric equations, as described for example in \cite{dworkpadic} and \cite{kloostermanFermat}, and reviewed in \cite{hypergeometric}.

\begin{proposition}\label{P:F4}
The primitive middle-dimensional cohomology group $H^2_{\textup{prim}}(X_{\Fsf_4, \psi},\C)$ has $21$ periods whose Picard--Fuchs equations are hypergeometric differential equations as follows:
\begin{itemize}
\item $3$ periods are annihilated by $D\bigl(\tfrac{1}{4}, \tfrac{1}{2}, \tfrac{3}{4} ; 1, 1, 1 \mid\psi^{-4} \bigr)$
\item $6$ periods are annihilated by $D\bigl(\tfrac{1}{4}, \tfrac{3}{4}; 1, \tfrac{1}{2} \mid\psi^{-4}\bigr)$
\item $12$ periods are annihilated by $D\bigl(\tfrac{1}{2};1 \mid\psi^{-4}\bigr)$.
\end{itemize}
\end{proposition}

Several of the periods for the $\Lsf_2\Lsf_2$ family satisfy trivial Picard--Fuchs equations. The rest are hypergeometric:

\begin{proposition}[\cite{hypergeometric}]\label{P:L2L2}
The group $H^2_{\textup{prim}}(X_{\Lsf_2\Lsf_2, \psi}, \C)$ has periods whose Picard--Fuchs equations are hypergeometric differential equations as follows:
\begin{itemize}
\item $3$ periods are annihilated by $D(\tfrac{1}{4}, \tfrac{1}{2}, \tfrac{3}{4} ; 1, 1, 1 \mid\psi^{-4})$
\item $2$ periods are annihilated by $D(\tfrac{1}{4}, \tfrac{3}{4} ; 1,  \tfrac12 \mid\psi^{4})$
\item$8$ periods are annihilated by $D(\tfrac{1}{8}, \tfrac{3}{8}, \tfrac58, \tfrac78 ; 0, \tfrac{1}{4}, \tfrac12, \tfrac34 \mid\psi^{4})$
\end{itemize}
\end{proposition}

To organize the corresponding point-counting information, we use an incomplete $L$-series.  Let $S$ be the set of bad primes for $X_\psi$, and define
\begin{equation}
L_S(X_{\diamond,\psi},s) = \prod_{p \not \in S} P_{X_\psi,p}(p^{-s})^{-1}
\end{equation}
which is convergent for $s \in \C$ in a right half-plane.

The main theorem of \cite{hypergeometric} writes $L_S(X_{\diamond,\psi},s)$ explicitly in terms of \emph{finite field hypergeometric functions} whose parameters are consistent with the Picard--Fuchs hypergeometric parameters.  Such functions have been studied by many authors, under different hypotheses on the parameters; \cite{hypergeometric} allows for a further weakening of the hypotheses.

Recall that one may define the rising factorials used in hypergeometric functions as ratios of gamma functions. The notion analogous to the gamma function over a finite field is a \emph{Gauss sum}.  Let $p$ be prime, let $q=p^r$, and abbreviate $q-1$ as $q^\times$.

\begin{definition}
Let $\omega \colon \F_q^{\times} \to \C^\times$ be a generator of the character group on $\F_q^{\times}$, and let $\Theta \colon \F_q \to \C^\times$ be a nontrivial additive character.  For $m \in \Z$, we define the Gauss sum $g(m)$ as
\begin{equation} 
g(m)= \sum_{x \in \F_q^{\times}} \omega(x)^m \Theta(x).
\end{equation}
\end{definition}

\noindent Here, the combination of the multiplicative and additive character is reminiscent of the definition of $\Gamma(z)$ using the integral $\int_0^\infty x^{z-1}e^{-x}\,dx$.

Let $\pmb{\alpha}=\{\alpha_1,\dots,\alpha_d\}$ and $\pmb{\beta}=\{\beta_1,\dots,\beta_d\}$ be multisets of $d$ rational numbers.  Suppose that $\pmb{\alpha}$ and $\pmb{\beta}$ are \emph{disjoint modulo $\Z$}, that is, $\alpha_i-\beta_j \not\in \Z$ for all $i,j \in \{1,\dots,d\}$. 

We now define a finite field hypergeometric sum. We follow work of Greene \cite{Greene} and Katz \cite[p.\ 258]{Katz} but normalize using the convention of \cite[Definition 3.2]{McCarthy} and Beukers--Cohen--Mellit \cite[Definition 1.1]{BCM}.

\begin{definition}\label{D:classicff} 
Suppose that
\begin{equation} \label{E:qqalpha}
\qq \alpha_i, \qq \beta_i \in \Z
\end{equation}
for all $i=1,\dots,d$.  For $t \in \F_q^\times$, we define a finite field hypergeometric sum by
\begin{equation} 
H_q(\pmb{\alpha}, \pmb{\beta} \mid t) = -\frac{1}{\qq} \sum_{m=0}^{q-2} \omega((-1)^dt)^m G(m+\pmb{\alpha}\qq,-m-\pmb{\beta}\qq)
\end{equation}
where
\begin{equation} \label{E:gmalphabeta}
G(m+\pmb{\alpha}\qq,-m-\pmb{\beta}\qq) = \prod_{i=1}^d \frac{g(m+ \alpha_i\qq)g(-m - \beta_i\qq)}{g(\alpha_i \qq)g(-\beta_i \qq)}
\end{equation}
for $m \in \Z$.
\end{definition}

The divisibility condition given in Equation~\ref{E:qqalpha} is restrictive. In \cite{BCM}, Beukers, Cohen, and Mellit gave an alternative definition of a finite field hypergeometric sum. We first define the \emph{field of definition} associated to hypergeometric parameters.

\begin{definition}[\cite{BCM}] \label{D:fieldofdef}
The field of definition $K_{\pmb{\alpha},\pmb{\beta}} \subset \C$ associated to $\pmb{\alpha},\pmb{\beta}$ is the field generated by the coefficients of the polynomials
\begin{equation}
\prod_{j=1}^d (x-e^{2\pi i \alpha_j}) \text{ and } \prod_{j=1}^d (x-e^{2\pi i \beta_j}). 
\end{equation}
In particular, if $\prod_{j=1}^d (x-e^{2\pi i \alpha_j})$ and $\prod_{j=1}^d (x-e^{2\pi i \beta_j})$ are products of cyclotomic polynomials, and thus have coefficients in $\mathbb{Z}$, we say $\pmb{\alpha},\pmb{\beta}$ is \emph{defined over $\mathbb{Q}$}.
\end{definition}

We say that $q$ is \emph{good} for $\pmb{\alpha},\pmb{\beta}$ if $q$ is coprime to the least common denominator of $\pmb{\alpha} \cup \pmb{\beta}$.

\begin{definition}[\cite{BCM}]\label{D:ffdefoverQ}
Suppose that $\pmb{\alpha},\pmb{\beta}$ are defined over $\Q$ and $q$ is good for $\pmb{\alpha},\pmb{\beta}$.  Choose $p_1, \ldots, p_r,q_1, \ldots, q_s \in \Z_{\geq 1}$ such that
\begin{equation}
\prod_{j=1}^d \frac{ (x-e^{2\pi i \alpha_j})}{(x-e^{2\pi i \beta_j})} = \frac{\prod_{j=1}^r x^{p_j} - 1}{\prod_{j=1}^s x^{q_j} - 1}.
\end{equation}
Let  $D(x)  =  \gcd(\prod_{j=1}^r (x^{p_j} - 1), \prod_{j=1}^s (x^{q_j} - 1))$, let $M = \bigl(\prod_{j=1}^r p_j^{p_j}\bigr) \bigl(\prod_{j=1}^s q_j^{-q_j}\bigr)$, set $\epsilon = (-1)^{\sum_{j=1}^s q_j}$, and let $s(m) \in \Z_{\geq 0}$ be the multiplicity of the root $e^{2\pi i m / \qq}$ in $D(x)$.  Finally, abbreviate
\begin{equation} 
g(\pmb{p}m,-\pmb{q}m) = g(p_1m) \cdots g(p_rm) g(-q_1m) \cdots g(-q_sm).
\end{equation}
For $t \in \F_q^\times$, define the finite field hypergeometric sum associated to $\pmb{\alpha}$ and $\pmb{\beta}$ by
\begin{equation} 
H_q(\pmb{\alpha}, \pmb{\beta} \mid t) = \frac{(-1)^{r+s}}{1-q} \sum_{m=0}^{q-2} q^{-s(0) + s(m)} g(\pmb{p}m,-\pmb{q}m)\omega(\epsilon M^{-1}t)^m.
\end{equation}
\end{definition}

By \cite[Theorem 1.3]{BCM}, Definitions~\ref{D:classicff} and \ref{D:ffdefoverQ} yield the same result in the cases where they both apply. However, situations may arise where neither definition applies directly. For example, analysis of the $\Fsf_1\Lsf_3$ family involves the hypergeometric parameters $\pmb{\alpha} = \{\frac{1}{14}, \frac{9}{14}, \frac{11}{14}\}$ and $\pmb{\beta} = \{0, \frac14, \frac34\}$.  We cannot use Definition~\ref{D:ffdefoverQ} since $(x-e^{2\pi i/14})(x-e^{18\pi i/14})(x-e^{22\pi i/14}) \not \in \Q[x]$.  When $q \equiv 1 \pmod{28}$, Definition~\ref{D:ffdefoverQ} applies; otherwise, it does not.  However, one may decompose these parameters into multisets where each of these definitions applies.

\begin{definition}[\cite{hypergeometric}]
We say that $q$ is \emph{splittable} for $\pmb{\alpha},\pmb{\beta}$ if there exist partitions 
\begin{equation} 
\pmb{\alpha} = \pmb{\alpha}_0 \sqcup \pmb{\alpha}' \textup{ and } \pmb{\beta} = \pmb{\beta}_0 \sqcup \pmb{\beta}' 
\end{equation}
where $\pmb{\alpha}_0,\pmb{\beta}_0$ are defined over $\Q$ and 
\[ \qq \alpha_i',\qq \beta_j' \in \mathbf{Z} \]
for all $\alpha_i' \in \pmb{\alpha}'$ and all $\beta_j' \in \pmb{\beta}'$.
\end{definition}

\begin{example}
Let $\pmb{\alpha} = \{\frac{1}{14}, \frac{9}{14}, \frac{11}{14}\}$ and $\pmb{\beta} = \{0, \frac14, \frac34\}$, and let $q$ be odd with $q \equiv 1 \pmod{7}$. Then $q$ is splittable for $\pmb{\alpha},\pmb{\beta}$: we may take $\pmb{\alpha}_0=\emptyset$, $\pmb{\alpha}'=\pmb{\alpha}$ and $\pmb{\beta}_0=\pmb{\beta}$, $\pmb{\beta}'=\emptyset$.
\end{example}

A hybrid definition of a finite field hypergeometric sum that applies in the splittable case is given in \cite{hypergeometric}.

With these definitions in hand, we may define an exponential generating series associated to a finite field hypergeometric series.

\begin{definition}
Let $t \in \F_q$. The finite field hypergeometric $L$-function associated to hypergeometric parameters $\pmb{\alpha}$ and $\pmb{\beta}$ is as follows.
\begin{equation}
L_\frakp(H_\frakp(\pmb{\alpha};\pmb{\beta}\mid t), T) = \exp\left(-\sum_{r=1}^{\infty} H_{\frakp^r}(\pmb{\alpha};\pmb{\beta}\mid t) \frac{T^r}{r} \right) \in K[[T]].
\end{equation}
\end{definition}

One may show that $L_\frakp(H_\frakp(\pmb{\alpha};\pmb{\beta}\mid t), T)$ is a polynomial. Furthermore, the degree of this polynomial matches the order of the hypergeometric differential equation with the same parameters.
 
The main theorem of \cite{hypergeometric} describes the $L$-functions for each of the five families in terms of finite field hypergeometric $L$-functions. We state the results for $\Fsf_4$ and $\Lsf_2 \Lsf_2$.

\begin{theorem}[\cite{hypergeometric}]\label{T:hypermainthm}
Let $t=\psi^{-4}$.
\begin{itemize}

\item For the Fermat pencil $\Fsf_4$,
\begin{align*}
L_S(X_{\Fsf_4,\psi}, s) &= L_S( H(\tfrac{1}{4}, \tfrac{1}{2}, \tfrac{3}{4}; 0, 0, 0\mid t), s) \\
&\qquad \cdot L_S( H(\tfrac{1}{4}, \tfrac{3}{4}; 0, \tfrac{1}{2} \mid t), s-1, \phi_{-1})^3 \\
&\qquad \cdot L_S( H(\tfrac{1}{2}; 0 \mid t) , \Q(\sqrt{-1}), s-1, \phi_{\sqrt{-1}})^6 
\end{align*}
where
\begin{align*}
\phi_{-1}(p) &=\legen{-1}{p} = (-1)^{(p-1)/2} &  & \text{ is associated to $\Q(\sqrt{-1}) \mid \Q$, and} \\
\phi_{\sqrt{-1}}(\frakp)&=\legen{i}{\frakp}=(-1)^{(\Nm(\frakp)-1)/4} & & \text{ is associated to $\Q(\zeta_8)\mid\Q(\sqrt{-1})$.}
\end{align*}

\item For the pencil $\Lsf_2 \Lsf_2$,
\begin{align*}
L_S(X_{\Lsf_2\Lsf_2,\psi}, s) &= L_S( H(\tfrac{1}{4}, \tfrac{1}{2}, \tfrac{3}{4}; 0, 0, 0\mid t), s) \\
&\qquad \cdot \zeta_{\Q(\sqrt{-1})}(s-1)^4  L_S( H(\tfrac{1}{4}, \tfrac{3}{4}; 0, \tfrac{1}{2} \mid t), s-1, \phi_{-1}) \\
&\qquad \cdot     L_S( H(\tfrac{1}{8}, \tfrac{3}{8}, \tfrac{5}{8}, \tfrac{7}{8}; 0, \tfrac{1}{4}, \tfrac{1}{2}, \tfrac{3}{4}\mid t), \Q(\sqrt{-1}), s-1, \phi_{\sqrt{-1}}\phi_{\psi})
\end{align*}
where
\[
\begin{aligned}
\phi_{\psi}(p) & = \legen{\psi}{p} & & \text{ is associated to $\Q(\sqrt{\psi})\mid\Q$}.
\end{aligned}
\]
\end{itemize}
\end{theorem}

Note that the finite field hypergeometric parameters in Theorem~\ref{T:hypermainthm} match the Picard--Fuchs hypergeometric parameters in Propositions~\ref{P:F4} and \ref{P:L2L2}. The trivial Picard--Fuchs equations for $\Lsf_2 \Lsf_2$ correspond to the zeta function $\zeta_{\Q(\sqrt{-1})}(s-1)^4 $. The proof of Theorem~\ref{T:hypermainthm} proceeds by explicit computation with finite field hypergeometric sums, using the corresponding Picard--Fuchs parameters as a guide.

Intuitively, one expects that as long as the discrete group of symmetries $SL(F_A)/J(F_A)$ commutes with the action of Frobenius, we will obtain subspaces of cohomology, with each block corresponding to both a Picard--Fuchs differential equation and a factor of $P_X$. However, such a factorization is only guaranteed over $\overline{\mathbb{Q}}$. To predict factors of $P_X$ over $\mathbb{Z}$, one must study the structure of the associated $L$-functions more carefully.

\begin{example}[\cite{hypergeometric}]
Let $Q_{\diamond,\psi,q}(T) = P_{X_{\diamond, \psi}}/R_\psi$. For smooth members of the $\Fsf_4$ and $\Lsf_2\Lsf_2$ families, the polynomials $Q_{\diamond,\psi,q}(T)$ factor over $\Q[T]$ as follows:
\begin{equation} \label{table:factdegrees}
\centering
\begin{tabular}[c]{c|c|c}
\textup{Family} & \textup{Factorization} & \textup{Hypothesis}  \\
\hline\hline
\rule{0pt}{2.5ex} \multirow{2}{*}{$\Fsf_4$} & $(\deg\, 2)^{3}(\deg\, 1)^{12}$ & $q\equiv 1 \psmod{4}$  \\
& $(\deg\, 2)^{3}(\deg\, 2)^6$ & $q \equiv 3 \psmod 4$ \\  \hline

\rule{0pt}{2.5ex} \multirow{2}{*}{$\Lsf_2\Lsf_2$} & $(1-qT)^8(\deg\, 2)(\deg\, 4)^2 $ & $q\equiv 1 \psmod{4}$ \\ 
 & $(1-q^2T^2)^4(\deg\, 2)(\deg\, 8) $ & $q\equiv 3 \psmod{4}$ \\ 
\end{tabular}
\end{equation}

\end{example}

The polynomials $Q_{\diamond,\psi,q}(T)$ may factor further, depending on specific values of $\psi$.

\bibliographystyle{amsalpha}

\end{document}